\documentclass[12pt]{article}
\usepackage{amsmath, amssymb, theorem}

\DeclareFontFamily{OT1}{rsfs}{} \DeclareFontShape{OT1}{rsfs}{m}{n}{
<-7> rsfs5 <7-10> rsfs7 <10-> rsfs10}{}
\DeclareMathAlphabet\mathcurl{OT1}{rsfs}{m}{n}

\theoremheaderfont\scshape
\newtheorem{theorem}{Theorem}[section]
\newtheorem{lemma}[theorem]{Lemma}
\newtheorem{corollary}[theorem]{Corollary}
\newtheorem{definition}[theorem]{Definition}
\newtheorem{proposition}[theorem]{Proposition}
\newtheorem{example}[theorem]{Example}

\theorembodyfont{\rmfamily}
\newtheorem{remark}[theorem]{Remark}

\newtheorem{assumption}[theorem]{Assumption}

\renewcommand\thefootnote{\fnsymbol{footnote}}

\newcommand\A{\mathcal A}
\newcommand\Rad{\mathcal R}
\newcommand\R{\mathbb R}
\newcommand\N{\mathbb N}
\newcommand\F{\mathcurl F}
\newcommand\E[2][\P]{\mathbb E_{#1}\left[ #2\right]}
\renewcommand\P{\mathbb P}
\newcommand\Q{\mathbb Q}

\renewcommand\d{\mathrm d}
\newcommand\qed{\hfill$\Box$}
\newcommand\ip[2]{\langle #1,#2 \rangle}
\newcommand\adm{\textup{adm}}
\newcommand\suprep{\overline\pi} 
\newcommand\radnik[1][\mathbb Q]{\frac{\d #1}{\d\mathbb P}}
\newcommand\normcl[2][\Q]{\overline{#2}^{L^1(#1)}}
\newcommand\ind{1\hspace{-2.5mm}1}
\newcommand\cone{{\cal{W}}} 

\newcommand{\DL}{F}
\newcommand{\DLCW}{\widehat{F_\Phi}}
\newcommand{\DLC}{F_\Phi}

\newcommand\comment[1]{}

\newcommand{\nohyphens}{\hyphenpenalty=10000\exhyphenpenalty=10000\relax}%
\makeatletter
\renewcommand{\@seccntformat}[1]{\csname the#1\endcsname.\hspace{1em}}%
\renewcommand\section{\@startsection{section}{1}{\z@}%
                {-3.5ex \@plus -1ex \@minus -.2ex}%
                {2.3ex \@plus.2ex}%
                {\setcounter{equation}0\bfseries\nohyphens}}%
\makeatother

\newenvironment{proof}[1][]{\noindent\textit{Proof#1.} }{\vskip\baselineskip}

\begin{document}
{\centering {\renewcommand\thefootnote{}\large\bfseries Geometry of
polar wedges and super-replication prices in incomplete financial
markets\footnote{The authors gratefully acknowledge support from
EPSRC grant no. GR/S80202/01.}

}



\par
\vspace{1em} {\scshape Frank Oertel}\par\itshape Department of
Mathematics\par University College Cork\par \vspace{1em} {\scshape
Mark P. Owen}
\par\itshape Department of Actuarial Mathematics and Statistics\par
Maxwell Institute for Mathematical Sciences and
Heriot-Watt University
}

\renewcommand{\abstractname}{}
\begin{abstract}\renewcommand\thefootnote{}\footnotesize
Consider a financial market in which an agent trades with
utility-induced restrictions on wealth. By introducing a general
convex-analytic framework which includes the class of umbrella
wedges in certain Riesz spaces (cf. Definition \ref{def:umbrella}
and Remark \ref{thm:risk measures and umbrellas}) and faces of
convex sets, consisting of probability measures (cf. Lemma
\ref{thm:likefrit} and Proposition \ref{thm:symmetric duality for
faces}), together with a duality theory for polar wedges, we provide
a representation of the super-replication price of an unbounded (but
sufficiently integrable) contingent claim that can be dominated
approximately by a zero-financed terminal wealth as the the supremum
of its discounted expectation under pricing measures which appear as
faces of a given set of separating probability measures (cf.
Proposition \ref{thm:symmetric duality} and Corollary
\ref{thm:general superreplication result}).

As a first direct application of our general framework, we see that
the utility-based super-replication price of an unbounded contingent
claim is equal to the supremum of its discounted expectations under
pricing measures with finite entropy (cf. Corollary \ref{thm:two}).
Central to our proof is the representation of a wedge $C_\Phi$ of
utility-based super-replicable contingent claims as the polar wedge
of the set of finite entropy separating measures. $C_\Phi$ is shown
to be the closure, under a relevant weak topology, of the wedge of
all (sufficiently integrable) contingent claims that can be
dominated by a zero-financed terminal wealth. Our general approach
shows, that those terminal wealths need {\it not} necessarily stem
from {\it admissible} trading strategies only.

We investigate also the natural dual of this result, and show in
particular that the polar wedge of $C_\Phi$ is the wedge generated
by separating measures with {\it finite loss-entropy}. For an agent
whose utility function is unbounded from above, the set of pricing
measures with finite loss-entropy can be larger than the set of
pricing measures with finite entropy. Indeed, we prove that the
former set is the closure of the latter under a suitable weak
topology (cf. Theorem \ref{thm:weakclose}).

The full two-sided polarity we achieve between measures and
contingent claims yields an economic justification for the use of
the wedge $C_\Phi$: the utility-based restrictions which this wedge
imposes on terminal wealth arise only from the investor's
preferences to asymptotically large negative wealth.

Finally, an application of our results to the special case of
admissible trading reveals that under a suitable weak topology
$C_\Phi$ is the closure of the wedge of all a.\,s. \!{\it bounded}
contingent claims that can be dominated by a terminal wealth
originating from admissible trading strategies.

\end{abstract}

\section{Introduction}\label{sec:motivation}
\noindent Although this paper primarily is written for a
mathematical audience who need not have a detailed knowledge of
mathematical finance (including the related terminology of
stochastic analysis), we occasionally have to use some specific
terminology which cannot be explained in detail here, due to the
limitation of space. Therefore, we would like to refer the reader to
the introductory overview references \cite{BaxRen96, DelbScha06,
Rom04, Scha04} and the further references therein.
\comment{ \noindent Although this paper can be read by a
mathematical audience without a strong knowledge of mathematical
finance (and the related terminology of stochastic analysis), we
occasionally have to use some specific expressions which cannot be
explained in detail here, due to the limitation of space;
particularly in this section. Therefore, we would like to refer the
reader to the introductory overview articles ... Hunt and Kennedy,
Baxter and Rennie \cite{Scha04} and the further references therein.
}

Firstly, let us revisit the ideal, non-realistic case, namely the
case of a \textit{complete} (financial) market.
In a complete market, there is a \textit{unique}
arbitrage-free
price of a given derivative security or contingent claim.
Its payoff (i.\,e., its value) is modelled as a random variable $X
\in L^1(\P)$, where $\P$ denotes the original (statistical)
probability measure. Given a finite time horizon $T>0$ and a
non-random, constant risk-free rate of interest $r > 0$, the unique
arbitrage-free price is then given by the expected discounted
terminal value $\E[\Q]{\exp(-rT)X} = \exp(-rT)\,\E[\Q]{X}$ which is
computed with respect to the \textit{risk-neutral probability
measure} $\Q$ under which the underlying asset's expected return
equals $r$ (i.\,e., under which the discounted asset price behaves
like a martingale) and which is equivalent to $\P$.
Due to the $\Q$-martingale property (respectively the $\Q$-sigma
martingale property) of the discounted marketed security price, $\Q$
is called an \textit{equivalent martingale measure $($EMM$)$}
(respectively an \textit{equivalent sigma martingale measure}
$($ESMM$)$).

However, due to the incorporation of features like price jumps,
transaction costs and illiquidity, financial markets in general do
not allow a complete replication of payoffs by trading in marketed
securities. They are \textit{incomplete}. The classic no-arbitrage
theory of valuation in complete markets, based on the unique price
of a self-financing replicating portfolio, is not adequate for
non-replicable payoffs in incomplete markets, where a perfect risk
transfer is not possible. In incomplete markets, there can be many
ESMMs, and one may ask for additional boundary criteria to select a
specific candidate $\Q^\ast$ of all available ESMMs in this market.
The expectation $\E[\Q^\ast]{\exp(-rT)X}$ is then the chosen
extension of the risk-neutral price function.

More generally, in an incomplete market with finite time horizon
$T$, the discounted price process of $d$ risky assets altogether is
modelled as a (non-continuous) $\R^d$-valued semimartingale $S =
(S_t)_{t\in [0, T]}$. Not every contingent claim $X$ can be
replicated by a self-financing trading strategy. For such contingent
claims there exists a whole interval of arbitrage-free prices, as
opposed to the case of a complete market, where there exists a
unique replication price, as roughly sketched above. An upper bound
for this price interval is the super-replication price
\begin{equation}
\suprep(X) := \inf\{x \in \R :
\text{there is an admissible } H  \text{ such that } X \leq x +
G_T(H)\},
\end{equation}
 where $G_T(H)
:= \int_0^T H_u\d S_u$ denotes the discounted cumulative gain or
loss at the time horizon $T$. Recall that an $\R^d$-valued
predictable process $H$ is called an {\it admissible trading
strategy} (cf. \cite{DelbScha94, DelbScha01}) if $H$ is
$S$-integrable, and there exists a constant $c \geq 0$ such that for
all $t\in[0,T]$,
\begin{equation}
  \int_0^tH_u\d S_u\geq -c, \qquad\P-\textrm{a.\,s}.
\end{equation}
Let
\begin{equation*} M_1^\adm(\P) : =
\{\Q \ll \P : K^{\adm} \subseteq L^1(\Q) \text{ and } \E[\Q]{X} \leq
0\ \text{ for all } X \in K^{\adm}\},
\end{equation*}
be the set of all {\it separating measures}, where
$K^{\adm}:=\{G_T(H) : H \text{ is admissible}\}$ denotes the wedge
of all terminal wealths originating from zero-financed admissible
trading strategies (cf. Section \ref{sec:frame} for the definition
of a wedge).

Although we will not not make use of it in this paper, we would like
to recall now some well-known representations of $M_1^\adm(\P)$
which depend upon the regularity (boundedness) of the semimartingale
$S$. The final statement
was shown in the well-known article \cite{DelbScha98}:
\begin{remark}
\begin{enumerate}
\item If $S$ is bounded then $M_1^\adm(\P)$ is the set of all
$\P$-absolutely continuous probability measures $\Q$ such that $S$
is a $\Q$-martingale. \item If $S$ is locally bounded then
$M_1^\adm(\P)=M^a(\P)$, where $M^a(\P)$ denotes the set of all
$\P$-absolutely continuous probability measures $\Q$ such that $S$
is a $\Q$-local martingale. \item In the general case (where $S$ may
not be locally bounded), if $M_1^\adm(\P) \neq \emptyset$ then
$M_1^\adm(\P)$ is the closure, in the topology induced by the total
variation norm, of the set $M_\sigma^a(\P)$ of all $\P$-absolutely
continuous probability measures $\Q$ such that $S$ is a
$\Q$-sigma-martingale.
\end{enumerate}
\end{remark}
If $S$ is locally bounded and the set $M^e(\P)$ of all
$\P$-\textit{equivalent} probability measures $\Q$ such that $S$ is
a $\Q$-local martingale is non-empty,
it is well-known (cf. \cite{DelbScha95}) that for $X$ bounded from
below
\begin{equation}\label{eqn:upeqlow}
\suprep(X) = \sup_{\Q\in M^e(\P)}\E[\Q]{X}.
\end{equation}
In the case of a general semimartingale $S$, equation
\eqref{eqn:upeqlow} holds if one substitutes the set $M^e(\P)$ with
the set $M^e_{\sigma}(\P)$, provided $M^e_{\sigma}(\P) \neq
\emptyset$\footnote{The condition $M^e_{\sigma}(\P) \neq \emptyset$
is equivalent to the (NFLVR) property (cf. \cite[Theorem
3.4]{DelbScha01}).} (cf. \cite[Theorem 5.12]{DelbScha98}).

It is already useful at this point to extend the definition of the
super-replication price to allow terminal wealths from an
\textit{arbitrary} fixed wedge $K$. Several candidates for $K$ may
be appropriate including \textit{admissible} strategies,
\textit{acceptable} strategies (cf. \cite{HugoKram04}) or
\textit{permissible} strategies (cf. \cite{OwenZitc07}).

Let $X \in L^0$ and $K \subseteq L^0$ be an \textit{arbitrarily}
given wedge. Assume that $A(X;K) \not= \emptyset$, where
\begin{eqnarray*}
A(X;K) & : = & \{x\in\R:X \leq x + G \textup{ for some } G \in
K\}\\
& \,\,= & \{ x\in\R : X - x \in K - {L}_+^0 \},
\end{eqnarray*}
and consider
\begin{equation*}
\suprep(X;K):=\inf(A(X;K)).
\end{equation*}
Note that $\suprep(X;K^{\adm}) = \suprep(X)$ if $A(X;K^{\adm}) \not=
\emptyset$ and that $\suprep(X;K_2) \leq \suprep(X;K_1)$ if $K_1
\subseteq K_2$ and $A(X;K_1) \not= \emptyset$.

Let us now assume that $M_1^\adm(\P) \not= \emptyset$. If one is
interested in pricing the claim $X$ by using separating measures
from the set $M_1^\adm(\P)$, it is natural to assume that $X \in
L(M_1^\adm(\P)) : =
\bigcap\limits_{{\Q} \in M_1^\adm(\P)}L^1(\Q)$. Recall that for such
$X$ we always have
\[
-\infty < \sup_{\Q\in M_1^\adm(\P)}\E[\Q]{X} \leq \suprep(X) <
+\infty\,.
\]
Note also that by construction
\[
K^\adm \subseteq L(M_1^{\adm}(\P))\,.
\]
Moreover, an easy calculation shows that the set $L(M_1^{\adm}(\P))$
contains all contingent claims $X$ in $L^0$ which are bounded from
below and satisfy $A(X;K^{\adm}) \not= \emptyset$.
\comment{
\begin{proof}
Let $-r \leq X, r \geq 0$. Then
\[
\vert X \vert = \vert (X+r) - r \vert \leq (X + r) + r \leq (x +
G_T(H)) +  2r. 
\]
\qed
\end{proof}
\begin{remark}
If $X \in L(M_1^{\adm}(\P))$
then
\[
-\infty < \sup_{\Q\in M_1^\adm(\P)}\E[\Q]{X} \leq \suprep(X) <
+\infty.
\]
\end{remark}
\begin{proof}
Let $x \in \R$ such that $X - x \in K^\adm - {L}^0_+$ and $\Q\in
M_1^\adm(\P)$. Then
\[
-\infty < - \E[\Q]{\vert X \vert } \leq \E[\Q]{X} = \E[\Q]{X-x} + x
\leq 0 + x\,.
\]
\qed
\end{proof}
}

However, one can construct $Y \in L(M_1^\adm(\P))$ such that
\[
\sup_{\Q \in M_1^\adm(\P)}\E[\Q]Y < \suprep(Y),
\]
suggesting that the use of admissible trading strategies is
unsuitable for super-replication of \textit{unbounded} claims (cf.
\cite{BiagFrit04}); the wedge $K^\adm$ is ``too small'' for the
purpose of super-replicating such contingent claims.

Let $X \in L(M_1^\adm(\P))$. Then
\[
\suprep(X; K^{\adm}) = \suprep(X; s_{L(M_1^\adm(\P))}(K^{\adm})),
\]
where
\begin{eqnarray*}
s_{L(M_1^\adm(\P))}(K^{\adm}) & = & \left\{G\in L(M_1^\adm(\P)): G
\leq U \textup{ for some } U \in K^{\adm} \right\} \\
& = & \bigcap\limits_{\Q\in M_1^\adm(\P)}\big(K^{\adm} -
L^1_+(\Q)\big)
\end{eqnarray*}
denotes the \textit{smallest} wedge $C \subseteq L(M_1^\adm(\P))$
such that \mbox{$-L(M_1^\adm(\P))_+ \subseteq C$} and $K^{\adm}
\subseteq C$ (cf. Definition \ref{def:umbrella} and
\eqref{eqn:umbrella rep of M based suprep price}). Cones of this
type will play an important role in our analysis. Their use namely
provides a significant extension of Theorem 5 in \cite{BiagFrit04}
in that it allows us to prove the super-replication result for
utility functions which
are not bounded from above (cf. Corollary \ref{thm:general
superreplication result} and Corollary \ref{thm:two}).

The following natural question appears. Is it possible to find an
\textit{enlarged} wedge $C \in L(M_1^\adm(\P))$, satisfying
$s_{L(M_1^\adm(\P))}(K^{\adm}) \subseteq C$, and a suitable
non-empty subset $M$ of $M_1^\adm(\P)$ (which may depend on
$K^\adm$), such that
\begin{equation}\label{eqn:arbcone} \suprep(X;C)=\sup_{\Q\in
M}\E[\Q]X\,? \end{equation} Can we even remove the admissibility
assumption and carry forward equation \eqref{eqn:arbcone} to a
larger class of wedges $K \subseteq L^0$ (cf. Corollary
\ref{thm:general superreplication result})?

In the admissible case, a partial answer to the question was
provided in \cite{BiagFrit04}, where preferences of the investor
were incorporated in the construction of the enlarged wedge by means
of the convex conjugate $\Phi$ of their utility function $U$. This
wedge was defined as
\begin{equation}\label{eqn:Cdefforadm} C_{\Phi}^{\adm}:=\bigcap_{\Q\in
M_\Phi^\adm}\normcl{K^{\adm} - L_+^1(\Q)},
\end{equation}
where
\begin{equation*} M_\Phi^\adm : =M_\Phi^\adm (\P) := \left\{\Q\in M_1^\adm(\P): \Phi \circ \radnik \in
L^1(\P)\right\} \end{equation*} denotes the set of pricing measures
with finite entropy.

Under the assumptions that the utility function $U$ has the
condition of Reasonable Asymptotic Elasticity at $-\infty$ and at
$+\infty$ and is bounded from above\footnote{i.e., if $\Phi(0) =
\lim\limits_{x \to \infty}U(x) < \infty$}, and if $\{\Q\in
M_\Phi^\adm : \Q \sim \P\}\neq \emptyset$, Biagini and Frittelli
showed that
\begin{equation}\label{eqn:phieq}
\suprep(X;C_{\Phi}^{\adm})=\sup_{\Q\in M_\Phi^\adm}\E[\Q]{X}.
\end{equation}
Their results hinge on the observation that the wedge
$C_{\Phi}^{\adm}$ of terminal wealths with zero initial endowment
and the wedge generated by the set $M_\Phi^\adm$ of pricing measures
with finite entropy are polar to one another. However for an
unbounded utility function, this bipolar relation fails, and no
results were obtained. One of our main results shows exactly what
happens in this situation.

To this end, we work within a general framework -- a duality theory
for wedges in Riesz spaces -- which not only allows us to consider
terminal wealths arising from non-admissible trading strategies.
If one wishes to obtain a similar bipolar relation here, one needs
to enlarge the set of finite entropy pricing measures, by
considering the set of pricing measures with {\it finite
loss-entropy}.

By identifying suitable representations of relevant {\it polar
wedges} in the context of our financial market, we provide a
transparent and short proof of Corollary \ref{thm:two}, which
generalises equation \eqref{eqn:phieq} for a utility function which
is unbounded from above. We are also able to relax a further
assumption of \cite{BiagFrit04}, namely that $\{\Q\in
M_\Phi:\Q\sim\P\}\neq \emptyset$.

\section{Preliminaries and Notations}\label{sec:frame}
In this section, we introduce some basic notation and terminology
which we will use throughout the paper. The scalar field for vector
spaces is assumed to be the real field $\R$ only, and most of our
notations and definitions from probability theory, convex analysis
and functional analysis are standard. We refer the reader to the
monographs \cite{Heus82, Jarc81, Wong93} for the necessary
background in functional analysis, and recommend the monographs
\cite{BorLew06, Rock72} for the basics of convex analysis.

From \cite{AlipTourky07}, let us recall that a {\it wedge} in a real
vector space $E$ is a non-empty convex set $C\subseteq E$ satisfying
$\lambda C\subseteq C$ for all $\lambda\geq 0$\footnote{In the
literature, a wedge quite often is denoted as a \textit{convex
cone}.}. Obviously, a non-empty subset $C$ of $E$ is a wedge if and
only if it is closed under addition and non-negative scalar
multiplication.
Let $S$ be an arbitrary non-empty subset of the vector space $E$.
Let $\cone(S)$ denote the smallest wedge in $E$ which contains $S$,
i.e., the wedge generated by $S$ (cf. \cite{AlipTourky07}). It is
easy to show that
\begin{equation}\label{eqn:cone_rep}
\cone(S) = \{\lambda x : \lambda \geq 0, x \in
\operatorname{co}(S)\} = \cone(\operatorname{co}(S)),
\end{equation}
where $\operatorname{co}(S)$ denotes the {\it{convex hull}} of $S$.
We therefore arrive at the following description of $[S]$,\, the
linear span of $S$, which we shall use in the proof of Proposition
\ref{thm:newinput}:
\begin{equation}\label{eqn:lin_hull}
[S] = \cone(S) - \cone(S).
\end{equation}
\begin{lemma}\label{thm:cone_of cone cap cx set}
Let $E$ be a real vector space, $S \subseteq E$ be convex and $T
\subseteq E$ be a wedge. If $S \cap T$ is non-empty then
\begin{equation}
\cone(S \cap T) = \cone(S) \cap T.
\end{equation}
\end{lemma}
\begin{proof}
Since $\cone(S) \cap T$ is a wedge which contains the non-empty set
$S \cap T$, it follows immediately that
\[
\cone(S \cap T) \subseteq \cone(S) \cap T\,.
\]
Therefore, without loss of generality, we may assume that $\cone(S)
\cap T \not= \{0\}$. So, let $x \in (\cone(S) \cap
T)\setminus\{0\}$. Since $0 \not= x \in \cone(S)$,
\eqref{eqn:cone_rep} implies that $x = \lambda w$, where $\lambda
> 0$ and $w \in \operatorname{co}(S) = S$. Since $x \in T$ and $T$ is a wedge,
it consequently follows that $w = \lambda^{-1}x \in S \cap T
\subseteq \cone(S \cap T)$. Hence, $(\cone(S) \cap T)\setminus\{0\}
\subseteq \cone(S \cap T)$, and the remaining inclusion follows.\qed
\end{proof}
In order to embed utility-based super-replication prices in a
mathematically concise framework, it is very useful to work in
vector lattices (or Riesz spaces). Standard references are for
instance \cite{AlipBurk06, LuxZaan71, Mey-Nie91, Schaef74}.
\begin{definition}\label{def:umbrella}
Let $(E, \leq)$ be a vector lattice and $C \subseteq E$ a wedge in
$E$. $C$ is an umbrella $($wedge$)$ in $E$ if $-E_{+} = \{x \in E :
x \leq 0\} \subseteq C$.
\end{definition}
Given an arbitrary wedge $K \subseteq E$, we denote by $s_E(K)$ the
umbrella hull, i.\,e., the smallest umbrella in $E$ which contains
$K$. Note that
\[
s_E(K) = \big\{x \in E : \exists \; g \in K \mbox{ s.\,t. } x \leq g
\big\} = K - E_+.
\]
Consequently, a wedge $C \subseteq E$ is an umbrella in $E$ if and
only if $C = s_E(C) = C - E_+$.

\begin{remark}\label{thm:risk measures and umbrellas}
Here, it should be noted that our umbrella wedges appear naturally
in another vibrant research topic in mathematical finance, which
also is based on convex analytic and functional analytic techniques,
namely in relation to (coherent) \textit{risk measures}. For
example, consider risk measures which are defined on the vector
lattice $E : = L^2(\Omega, \F, \P)$ (cf. \cite{RUM02}).
Reinterpreting the conditions (A1), (A2), (A3) and (A5) of Theorem 2
in \cite{RUM02} and recalling Definition \ref{def:umbrella}, we
immediately obtain
\begin{proposition}
\textit{Let $E : = L^2(\Omega, \F, \P)$ and $\A$ be an arbitrary
non-empty subset of $E$. Then the following statements are
equivalent:
\begin{itemize}
\item[$(i)$] There exists a lower semicontinuous coherent risk
measure $\rho : E \longrightarrow \R \cup \{+\infty\}$, such that
\[
\A = \{X \in E : \rho(X) \leq 0\}.
\]
\item[$(ii)$] $-\A = \{-X : X \in \A\}$ is a closed umbrella wedge in
E.
\end{itemize}}
\end{proposition}
\end{remark}
\begin{proposition}\label{thm:rep of umbrella cones in sublattices}
Let $(E_\alpha)_{\alpha \in \mbox{A}}$ be a family of vector
sublattices of a vector lattice $L$. Then
\comment{Definition: A \textit{sublattice} $M$ of a lattice L is a
non-empty subset $M$ of $L$ which is a lattice with the same meet
and join operations as $L$. That is, if $L$ is a lattice and $M$ is
a non-empty subset of $L$ such that for every pair of elements $a,
b$ in $M$ both $a \wedge b$ and $a \vee b$ are in $M$, then
M is a sublattice of L.\\
\\
A sublattice $M$ of a lattice $L$ is a \textit{convex sublattice} of
$L$, if $x \leq  z \leq y$ and $x, y \in M$ implies that $z$ belongs
to M, for all elements $x, y, z$ in $L$.}
$E : = \bigcap_{\alpha \in \mbox{A}} E_\alpha$ is a vector
sublattice of $L$. Let $K \subseteq E$ be an arbitrary wedge in $E$.
Then
\[
s_{E}(K) = \bigcap_{\alpha \in \mbox{A}}s_{E_\alpha}(K) =
\bigcap_{\alpha \in \mbox{A}}\big(K - (E_\alpha)_+\big)\,.
\]
\end{proposition}
\begin{proof}
Let $C$ be an arbitrary umbrella wedge in $E$ such that $K \subseteq
C$. Then
\[
\bigcap\limits_{\alpha \in \mbox{A}}s_{E_\alpha}(K) \subseteq
\bigcap\limits_{\alpha \in \mbox{A}}\big(s_{E_\alpha}(K) \cap E\big)
\subseteq \bigcap\limits_{\alpha \in \mbox{A}}\big(s_{E_\alpha}(C)
\cap E\big) \subseteq \bigcap\limits_{\alpha \in \mbox{A}}s_{E}(C) =
C\,.
\]
Hence, $\bigcap\limits_{\alpha \in \mbox{A}}s_{E_\alpha}(K)
\subseteq s_{E}(K)$. The other inclusion is trivial. \qed
\end{proof}
\begin{example}
Assume that $M_1^\adm(\P) \not= \emptyset$. Set
$E : = \bigcap_{\Q \in M_1^\adm(\P)}L^1(\Q)$. Then $K^\adm \subseteq
E$, and
\[
s_{E}(K^\adm) = \bigcap_{\Q \in M_1^\adm(\P)}\big(K^\adm -
L^1_+(\Q)\big)\,.
\]
\end{example}
For the convenience of the reader we recall some facts of basic
duality theory, for which we adopt the approach of \cite{Heus82}.
Let $E$ and $\DL$ be vector spaces over $\R$. If there exists a
bilinear form $\ip\cdot\cdot$ on $E \times \DL$, the pair of vector
spaces $(E, \DL)$ is called a {\it{bilinear system}} (with respect
to the bilinear form $\ip\cdot\cdot$). A bilinear system $(E, \DL)$
is called a {\it{left dual system}} if $\ip z{w} = 0$ for all $z \in
E$ implies that $w = 0$, and $(E, \DL)$ is called a {\it{right dual
system}} if $\ip z{w} = 0$ for all $w \in E$ implies that $z = 0$.
If $(E, \DL)$ is a left dual system, then $\DL$ is the topological
dual of $E$ under the weak topology $\sigma(E,\DL)$ (see
\cite[\S{}82]{Heus82}):
\begin{equation*}
\DL = (E, \sigma(E,\DL))'.
\end{equation*}
The topology $\sigma(E,\DL)$ is Hausdorff if and only if $(E, \DL)$
is a right dual system. A bilinear system which is both a left dual
system and a right dual system is called a {\it dual system} (or a
{\it duality}). If $(E, \DL)$ is a dual system, then each of the
vector spaces $E, \DL$ is the topological dual of the other with
respect to the weak topologies, which are both then Hausdorff.

Let $(E, \DL)$ be an arbitrary bilinear system. For a non-empty set
$A\subseteq E$ we define its {\it polar wedge} $A^\lhd \subseteq
\DL$ by
\begin{equation*}
A^\lhd:=\{w\in \DL:\ip z{w} \le 0 \quad\forall z\in A\}.
\end{equation*}
For a non-empty set $B\subseteq \DL$ we define, in a similar way,
its polar wedge $B^\lhd \subseteq E$ by
\begin{equation*}
B^\lhd:=\{z \in E:\ip z{w} \le 0 \quad\forall w \in B\}.
\end{equation*}
Clearly, $A^\lhd$ is a wedge and $A^\lhd \subseteq
A^{\displaystyle\circ}$, where $A^{\displaystyle\circ} := \{w\in
\DL:\ip z{w} \le 1 \quad \forall z\in A\}$ denotes the {\it{polar}}
of $A$ (cf. \cite[Section 0.7]{Wong93}).
\begin{remark} Note that the definition of a polar is not
handled univocally! Some authors prefer to define the polar of $A$
as the set $\{w\in \DL: | \ip z{w} | \leq 1 \quad\forall z\in A\}$
(cf. e.\,g. \cite{Heus82}). However, if the set $A$ is circled,
i.\,e., if $\{\lambda x : \lambda \in [-1, 1], x \in A\} \subseteq
A$, both definitions coincide. Throughout the article, we will use
the definition above including the related version of the bipolar
theorem (see \cite[Theorem 0.8]{Wong93}).
\end{remark}
If in addition $A$ is a wedge, it easily follows that $A^\lhd =
A^{\displaystyle\circ}$. The next result shows why we call $A^\lhd$
a ``polar wedge'':
\begin{proposition}\label{thm: representation of a polar_cone}
Let $(E, \DL)$ be an arbitrary bilinear system of real vector spaces
and $A$ be an arbitrary non-empty subset of $E$. Then:
\begin{equation*}
A^\lhd = (\cone(A))^{\displaystyle\circ}.
\end{equation*}
In particular, $A^\lhd$ is $\sigma(\DL,E)$-closed.
\end{proposition}
\begin{proof}
It suffices to prove that $A^\lhd$ is contained in
$(\cone(A))^{\displaystyle\lhd}=(\cone(A))^{\displaystyle\circ}$.
Let $w \in A^\lhd$ and $z \in \cone(A)$. Due to
\eqref{eqn:cone_rep}, $z$ can be written as $z = \lambda u$, where
$\lambda \geq 0$ and $u \in \operatorname{co}(A)$. Thus, there exist
$n \in \N$, $\lambda_{1} \geq 0, \ldots, \lambda_{n} \geq 0$ and
$u_1, \ldots, u_n \in A$ such that $\sum_{k =1}^{n}\lambda_{k} = 1$
and $z = \sum_{k =1}^{n}\lambda \lambda_{k} u_k$. Since $w \in
A^\lhd$, the claim follows.\qed
\end{proof}
Consequently, a direct application of the bipolar theorem (see
\cite[Theorem 0.8]{Wong93}) leads to the following result which we
shall use in the proof of Theorem \ref{thm:weakclose}:
\begin{proposition}\label{thm:bipolar_cones}
Let $(E, \DL)$ be an arbitrary bilinear system of real vector spaces
and $A$ be an arbitrary non-empty subset of $E$. Then:
\begin{equation*}
A^{\lhd\lhd} = \overline{\cone(A)}^{\sigma(E, \DL)}.
\end{equation*}
In other words, $A^{\lhd\lhd}$ is the smallest $\sigma(E,
\DL)$-closed wedge which contains $A$. If in addition $A$ is a
wedge, then $A^{\lhd\lhd} = \overline{A}^{\sigma(E, \DL)}$.
\end{proposition}
\begin{lemma}\label{thm:permanence properties of dual cones}
Let $(E, \DL)$ be an arbitrary bilinear system and
$(A_\gamma)_{\gamma \in \Gamma}$ be a family of non-empty subsets of
$E$. Then
\begin{itemize}
\item[$($i$)$] $\Big(\bigcup\limits_{\gamma \in \Gamma} A_\gamma\Big)^\lhd
= \bigcap\limits_{\gamma \in \Gamma} A_\gamma^\lhd$ \mbox{ and }
\item[$($ii$)$] $ \Big(\bigcap\limits_{\gamma \in \Gamma} \overline{\cone{(A_\gamma)}}^{\sigma(E,
\DL)}\Big)^\lhd = \overline{\cone{\Big(\bigcup\limits_{\gamma \in
\Gamma} A_\gamma^\lhd\Big)}}^{\sigma(E, \DL)}$.
\end{itemize}
\end{lemma}
\begin{proof}
Statement (i) follows immediately from the definition of a polar
wedge.

If we apply statement (i) to the family $(B_\gamma)_{\gamma \in
\Gamma}$, where $B_\gamma : = A_\gamma^\lhd$, Proposition
\ref{thm:bipolar_cones} implies that
\[
\Big(\bigcap\limits_{\gamma \in \Gamma}
\overline{\cone{(A_\gamma)}}^{\sigma(E, \DL)}\Big)^\lhd   =
\Big(\bigcap\limits_{\gamma \in \Gamma} B_\gamma^{\lhd}\Big)^\lhd
\stackrel{(i)}{=} \Big(\bigcup\limits_{\gamma \in \Gamma}
B_\gamma\Big)^{\lhd\lhd}.
\]
Again, due to Proposition \ref{thm:bipolar_cones}, we obtain
\[
\Big(\bigcup\limits_{\gamma \in \Gamma} B_\gamma\Big)^{\lhd\lhd} =
\overline{\cone{\Big(\bigcup\limits_{\gamma \in \Gamma}
B_\gamma\Big)}}^{\sigma(E, \DL)} =
\overline{\cone{\Big(\bigcup\limits_{\gamma \in \Gamma}
A_\gamma^\lhd\Big)}}^{\sigma(E, \DL)},
\]
and statement (ii) follows.\qed
\end{proof}
In the appendix the interested reader will find another very
interesting application of this result which shows that an
infinite-dimensional version of Farkas' Lemma is true if and only if
suitable linear images of positive wedges are weakly closed (cf.
Theorem \ref{thm:Farkas_infty}).

We conclude this preparatory section with a technical lemma which we
shall use several times.
\begin{lemma}\label{thm:truncation_inequality}
Let $(\Omega, \mathcurl{F}, \P)$
be a probability space and let $Y_0, Y$ and $Y_1$ be real random
variables on $\Omega$ such that $0 \leq Y_0 \leq Y \leq Y_1$
almost surely. Let $f : [0, \infty) \longrightarrow [0,
\infty]$ be a convex function. Then for every $a \in [0, \infty]$
\begin{equation*}
\E{f(Y)\,\ind_{\{Y < a\}}} \leq \E{f(Y_0)} +
\min\{\E{f(Y_1)},f(a)\}.
\end{equation*}
\end{lemma}

\begin{proof}
Let $a \in [0, \infty]$ and put
\begin{equation*}
\gamma^\ast := \inf\{y \geq 0 : f \textrm{ is increasing on }
(y, \infty) \}.
\end{equation*}
Obviously, we only have to consider the case $\gamma^\ast < \infty$.
In this case, the convex function $f$ is non-increasing on $[0,
\gamma^\ast)$ and increasing on $[\gamma^\ast, \infty)$. Moreover
\begin{equation*}
f(Y)\,\ind_{\{Y < a\}} = f(Y)\,\ind_{\{Y < a \wedge
\gamma^\ast\}} + f(Y)\,\ind_{\{\gamma^\ast \leq Y < a\}}\,,
\end{equation*}
and the claim follows. \qed
\end{proof}

\section{The Market Model}\label{sec:model}
We now describe in detail the market model. For the
necessary background in mathematical finance and stochastic
analysis, we refer the reader to the introductory
monographs \cite{Kleb98,Prot03} and the survey
article \cite{DelbScha01}.

Let $(\Omega,\mathcurl F,(\mathcurl F_t)_{t\in[0,T]},\P)$ be a
filtered probability space, in which the filtration satisfies the
usual conditions of right continuity and completeness, and the
time horizon $T$ is assumed to be finite. 

We model the discounted price process of $d$ risky assets as an
$\R^d$-valued semimartingale $S=(S_t)_{t\in[0,T]}$. Let $K$ be an
{\it arbitrary} wedge in $L^0$ consisting of attainable terminal
wealths. Note that we {\it are not requiring} $K$ to be the wedge of
those attainable terminal wealths which arise from {\it admissible}
trading strategies.

We shall be particularly interested in super-replication for an
investor whose preferences are expressed via a utility function.
When considering the permissible trading strategies, it is important
to take into consideration the investor's wealth preferences. We
assume that the investor has a utility function
$U:(a,\infty)\rightarrow\R$, where
$a\in[-\infty,\infty)$\footnote{We allow $a$ to take any value in
$[-\infty,\infty)= \{-\infty\} \cup \R$, but we shall be most
interested in the case where $a=-\infty$.}, which is strictly
increasing, strictly concave, (continuously) differentiable, and
satisfies the Inada conditions
\begin{equation}\label{eqn:inada}
\lim_{x\downarrow a}U'(x)=\infty,\qquad
\lim_{x\uparrow\infty}U'(x)=0.
\end{equation}
As usual, we assume that the utility function has reasonable
asymptotic elasticity (see \cite{Scha00}). We shall formulate this
assumption however in Section \ref{sec:sep}, in terms of a growth
condition on the slightly more general convex (conjugate) function
$\Phi$.

\section{The Separating Measures}\label{sec:sep}
Throughout the paper, $K$ always denotes a fixed but
\textit{arbitrarily chosen wedge in $L^0$}. Define
\begin{eqnarray*}
M_1 & : = & M_1(\P; K) \\
& := & \{\Q\ll\P: \Q \textup{ is a probability measure, }K \subseteq
L^1(\Q) \textup{ and } \E[\Q]{X} \leq 0 \textup{ for all }X\in K\}.
\end{eqnarray*}
By construction, it follows immediately that
\begin{equation}\label{eqn:cone K in lattice}
K \subseteq \bigcap_{\Q\in M_1}L^1(\Q) \subseteq \bigcap_{\Q\in
M}L^1(\Q) \end{equation} for any non-empty subset $M$ of $M_1$. Note
that $M_1$ consists of probability measures (a fact needed in the
proof of Proposition \ref{thm:newinput} and Lemma \ref{thm:rep of
dual umbrella cones}).
Consider the following two convex sets:
\[
\A(\P) : = \big\{\Q : \Q \mbox{ is a probability measure and } \Q
\ll \P \big\}
\]
and
\[
\Rad(\P) : = \big\{ X : X \in L^1_+(\P),\,\E{X} = 1 \big\} \subseteq
L^1(\P).
\]
Then
\begin{eqnarray*}\label{eqn:canonical identification}
R : & \A(\P)  \rightarrow \Rad(\P)\\
{} & \hspace{2mm} \Q \mapsto \displaystyle{\radnik[\Q]}
\end{eqnarray*}
is one-to-one surjective. The inverse mapping is given by
\begin{eqnarray*}\label{eqn:canonical identification inverse}
L : & \Rad(\P)  \rightarrow \A(\P),\\
{} & \hspace{2mm} X \mapsto \Q_X
\end{eqnarray*}
where $\Q_X(A) : = \E{\ind_A\,X} = \int_A X \mbox{d}\P$ for any $A
\in \F$. Thus, $M_1 \subseteq \A(\P)$ and
\begin{equation}\label{eqn:altdefM1}
R(M_1)  =  \big\{Y\in L_+^1(\P): \E Y=1, XY\in L^1(\P) \text{ and
}\E{XY}\le 0\text{ for all }X\in K \big\}.
\end{equation}

If it is not mentioned explicitly, then throughout the article, we
let \mbox{$\Phi:[0,\infty)\rightarrow(-\infty,\infty]$} denote an
arbitrary convex function, where $\Phi$ may take the value $\infty$
only at $0$.\footnote{In other words, $\Phi\big((0, \infty)\big)
\subseteq \R$.}. Note that we do not necessarily require $\Phi$ to
be the convex conjugate of any utility function.

The set of pricing measures with finite entropy is defined as
\begin{align}
  M_\Phi :=& \left\{\Q\in M_1:\Phi\left(\radnik\right)\in
  L^1(\P)\right\}\notag\\
  =& \left\{\Q\in M_1:\E{\Phi^+\left(\radnik\right)}<\infty\right\}.
\end{align}
By construction, $M_\Phi$ depends on the choice of the wedge $K$.

We now state an important growth condition on the convex function
$\Phi$, which is related to the condition of reasonable asymptotic
elasticity found in \cite{Scha00}.
\begin{assumption}\label{thm:growth}
Given any interval $[\lambda_0,\lambda_1]\subseteq(0,\infty)$
there exist $\alpha>0$ and $\beta>0$ such that
\begin{equation}\label{eqn:growth}
  \Phi^+(\lambda y)\le\alpha\Phi^+(y)+\beta(y+1),
\end{equation}
for all $y>0$ and all $\lambda\in[\lambda_0,\lambda_1]$.
\end{assumption}
Let $\Phi$ satisfy the growth condition \eqref{eqn:growth}. Then
obviously $\Phi^+(\lambda \radnik{}) \in L^1(\P)$ for any $\lambda
>0$ and any $\Q \in M_\Phi$.

Primarily, we are interested in the case where $\Phi$ is the convex
conjugate of a utility function $U$, i.\,e., the Legendre conjugate
of the convex function $-U(- \cdot): (-\infty, -a) \longrightarrow
\R$ (cf. \cite[\S 26]{Rock72}):
\begin{equation}\label{eqn:conj}
  \Phi(y):=\sup_{x>a}\{U(x)-xy\},\hspace{1cm}y \geq 0.
\end{equation}
In this case, the growth condition \eqref{eqn:growth} and the
Reasonable Asymptotic Elasticity of $U$, introduced by
Schachermayer, are equivalent (cf. \cite{FritGia02}).

It is well-known that under the Inada conditions \eqref{eqn:inada}
the convex conjugate $\Phi$ is strictly convex and continuously
differentiable on $(0, \infty)$, satisfying
\begin{equation}\label{eqn:asymptotic condition}
\lim\limits_{y \rightarrow \infty}\Phi(y) = \lim\limits_{x
\rightarrow a}U(x) \hspace{0.5cm} \text{and} \hspace{0.5cm} {\Phi}'
= - I,
\end{equation}
where $I :(0, \infty) \longrightarrow  (a, \infty)$ denotes the
inverse function of $U'$. Moreover,
\begin{equation*}
\Phi(y) = U(I(y)) - y I(y)
\end{equation*}
for all $y >0$ (cf. \cite[\S 26]{Rock72}).

Throughout the article we shall assume the following
\begin{assumption}\label{thm:utarbfr}
$M_\Phi\neq\emptyset$.
\end{assumption}
This assumption rules out ``utility-based arbitrage strategies''
(see \cite[Section 1.2]{BiagFrit04}), but not necessarily a free
lunch with vanishing risk. As our analysis shows, it is not
necessary to assume the stronger condition
\begin{equation}
  \{\Q\in M_\Phi:\Q\sim\P\}\neq\emptyset
\end{equation}
which was required in \cite{BiagFrit04}.

\begin{definition}
We say that a measure $\Q\ll\P$ has \textup{finite loss-entropy} if
\begin{equation}\label{eqn:lossdefn}
 \E{\Phi^+\left(\radnik\right)\ind_{\{\radnik\ge1\}}}<\infty.
\end{equation}
The set of pricing measures with finite loss-entropy is therefore
given by
\begin{equation}\label{eqn:measfinitele}
\widehat M_\Phi := \left\{\Q\in
M_1:\E{\Phi^+\left(\radnik\right)\ind_{\{\radnik\ge1\}}}<\infty\right\}.
\end{equation}
\end{definition}

Clearly, $M_\Phi \subseteq \widehat M_\Phi$. Since
\[
0 \leq \Phi^+(\radnik)\ind_{\{a \leq \radnik <
b\}}\leq\max\{\Phi^+(a), \Phi^+(b)\}<\infty
\]
for any $\Q \ll \P$ and $0<a<b$, the choice of the constant $1$ in
equation \eqref{eqn:lossdefn} is arbitrary; we could actually choose
any positive number, and the set in \eqref{eqn:measfinitele} would
not change. We use the terminology ``loss-entropy'' because for
events with large $\radnik$, an inspection of equation
\eqref{eqn:asymptotic condition} shows that
$\Phi\left(\radnik\right)$ is related to the value of the utility
function $U(x)$ where $x$ is close to the critical wealth $a$.
Typically, pricing measures $\Q\in M_1$ give large probabilities
(relative to the real world measure $\P$) to large negative asset
prices.

We now define a modification $\widehat\Phi$ of $\Phi$ which is
finite at $0$, but remains convex and satisfies
$\widehat\Phi\le\Phi$. Since the left-sided derivative
$l^\ast:=(D^-\Phi)(1)$ of the convex function $\Phi$ at $y=1$ always
exists, it follows that
\begin{equation}\label{eqn:defphihat}
\widehat\Phi(y):=\begin{cases}\Phi(y),&\text{if
  }y\ge1\\\Phi(1)+l^\ast(y-1)&\text{if }0 \leq y<1\end{cases}
\end{equation}
defines a function $\widehat\Phi$ on $[0, \infty)$ with values in
$\R$. If $\Phi$ satisfies the growth condition \eqref{eqn:growth},
so does $\widehat\Phi$. Moreover,
\begin{equation*}
{\widehat{M}}_\Phi = M_{\widehat\Phi}.
\end{equation*}
\begin{remark}\label{thm:pricingremarks}
\begin{enumerate}
\item Since $\Phi$ and $\widehat\Phi$ are convex, it follows easily
that $M_\Phi$ and $\widehat M_\Phi$ are convex. \item If
$\lim_{y\rightarrow\infty}(D^-\Phi)(y)<\infty$ (i.e., $\Phi$ is
asymptotically linear as $y\rightarrow\infty$), then $\widehat
M_\Phi=M_1$, because then $\widehat\Phi$ can be bounded above by a
linear function (see also Remark \ref{thm:laterremark}). \item If
$\Phi(0) < \infty$, then $\widehat M_\Phi = M_\Phi$, because in this
case $ 0 \leq \Phi^+(\radnik)\ind_{\{0 \leq \radnik < 1\}} \leq
\max\{\Phi^+(0), \Phi^+(1)\} < \infty $ for any $\Q \ll \P$.
\end{enumerate}
\end{remark}
Let us recall that a convex subset $F$ of a given convex set $C$ is
called a \textit{face of $C$} if $\alpha x + (1 - \alpha)y \in F$
with $x, y \in C$ and $0 <\alpha <1$ imply $x, y \in F$ (cf.
\cite{AlipTourky07}).
\begin{lemma}\label{thm:likefrit} If $\Phi$
satisfies the growth condition \eqref{eqn:growth} then the following
statements hold:
\begin{enumerate}
\item $\widehat M_\Phi$ is a face in $M_1$.
\item Let $\Q_0,\Q_1\in M_1$, $0 < x < 1$, and assume that $\Q : = x\Q_0+(1-x)\Q_1
\in\widehat M_\Phi$. If $\Q_0\in M_\Phi$ or $\Q_1\in M_\Phi$ then
$\Q\in M_\Phi$.
\end{enumerate}
\end{lemma}
\begin{proof}
Let $\Q_0,\Q_1\in M_1$, $0 < x < 1$ and define $\Q : =
x\Q_0+(1-x)\Q_1$. To prove statement (i), note that since $\widehat
M_\Phi$ is convex, we only have to show that $\Q\in\widehat M_\Phi$
implies that $\Q_0, \Q_1\in\widehat M_\Phi$. Suppose that
$\Q\in\widehat M_\Phi = M_{\widehat\Phi}$. Since $0 \leq
\radnik[\Q_0] = \frac{1}{x} \radnik - \frac{(1-x)}{x}\radnik[\Q_1]
\leq \frac{1}{x}\radnik \; \P-\textrm{a.\,s.}$, Lemma
\ref{thm:truncation_inequality} and the growth condition
\eqref{eqn:growth} of $\widehat\Phi$ imply that
\begin{equation*}
\E{\widehat\Phi^+\left(\radnik[\Q_0]\right)} \leq
\widehat\Phi^+(0)+\E{\widehat\Phi^+\left(\frac1x\radnik\right)}<\infty.
\end{equation*}
Hence, $\Q_0\in\widehat M_\Phi$. Similarly, $\Q_1\in\widehat
M_\Phi$.

To prove statement (ii), we may suppose without loss of generality
that $\Q_0\in M_\Phi$. Since $x \radnik[\Q_0] \leq \radnik$
$\P$-a.\,s., Lemma \ref{thm:truncation_inequality} and the growth
condition \eqref{eqn:growth} of $\Phi$ imply that
\begin{equation*} \E{\Phi^+\left(\radnik\right)\ind_{\{\radnik < 1\}}}
 \le \E{\Phi^+\left(x\radnik[\Q_0]\right)} + \Phi^+(1) <\infty. \end{equation*}
Since $\Q\in\widehat M_\Phi$, it therefore follows that $\Q\in
M_\Phi$.\qed
\end{proof}
We now include a general representation result which holds for
arbitrary nonempty faces of convex sets, consisting of probability
measures (such as it is the case for $M_1$). It is central to our
analysis.
\begin{proposition}\label{thm:faces of M1 and their linear hull}
Let $M$ be a face of $M_1$. If $M \not= \emptyset$ then
\[
[M] \cap M_1 = M\,.
\]
\end{proposition}
\begin{proof}
Clearly, $\emptyset \not= M \subseteq [M] \cap M_1$. Let $\Q \in [M]
\cap M_1$. First, due to \eqref{eqn:lin_hull}, it follows that
\begin{equation*}
[M] = \cone(M) - \cone(M).
\end{equation*}
Therefore, since $\Q$ in particular is a probability measure, $\Q$
can be written as $\Q=(1+\beta)\Q_1-\beta\Q_0$ where $\Q_0,\Q_1\in M
\subseteq M_1$ and $\beta\ge0$. Thus \mbox{$\Q_1 = \frac1{1+\beta}\Q
+ \frac\beta{1+\beta}\Q_0$}, and Lemma \ref{thm:likefrit} implies
that $\Q\in M$.\qed
\end{proof}
\begin{proposition}\label{thm:newinput} If $\Phi$ satisfies the
growth condition \eqref{eqn:growth} and $M_\Phi\neq\emptyset$ then
\begin{equation*}
\widehat M_\Phi = [M_\Phi] \cap M_1.
\end{equation*}
\end{proposition}
\begin{proof}
Firstly, due to Lemma \ref{thm:likefrit}, we have
\[
[M_\Phi] \cap M_1 \subseteq [\widehat M_\Phi] \cap M_1 = \widehat
M_\Phi.
\]
Let $\Q_1\in\widehat M_\Phi$. Choose any $\Q_0\in M_\Phi$ and define
$\Q : = \frac12\Q_0+\frac12\Q_1$. From Lemma \ref{thm:likefrit} we
see that $\Q\in M_\Phi$. Hence, $\Q_1=2\Q-\Q_0\in [M_\Phi]$.\qed
\end{proof}
\section{Main Results: A general Approach by Duality}\label{sec:results}
In the following, we extend the model approach of \cite{BiagFrit04},
which does not only allow us to substitute the wedge $K^\adm$ of
admissible strategies by more general candidates including
the wedge of acceptable strategies (cf. \cite{HugoKram04}), or the
wedge of permissible strategies (cf. \cite{OwenZitc07}). This
extension also enables to replace the set $M_\Phi$ by other suitable
non-empty subsets of $M_1$ (cf. Corollary \ref{thm:general
superreplication result}) and reveals the geometry of  a suitable
class of wedges which appears naturally in the analysis of
superreplication prices for unbounded contingent claims.

To this end, let $K$ be a fixed wedge in $L^0$ and $M$ be a fixed
non-empty subset of $M_1 \equiv M_1(\P; K) \subseteq \A(\P)$.
\comment{ By identifying probability measures $\Q\ll\P$ with their
Radon-Nikodym derivatives $\radnik\in L^1_+(\P)$ of norm $1$ (via
the mapping $R$, respectively $L$), we may \textit{identify} the
linear span $[M]$ with $\DL$, where $\DL$ denotes the subspace of
$L^1(\P)$, formed by taking the linear span of the Radon-Nikodym
derivatives of all probability measures in $M$. To this end, we only
have to consider the canonical linear extension of $R\vert_M$ which
is defined on $[M]$.}
We consider the pair of vector spaces %
\begin{equation}\label{eqn:bilinear system} E : = \bigcap_{\Q\in M}L^1(\Q)
\qquad\textup{and}\qquad \DL : = [R(M)] = \bigg[\radnik : \Q \in M
\bigg] \subseteq L^1(\P),
\end{equation}
where $R(\Q) = \displaystyle\radnik$ for all $\Q \in \A(\P)$ (cf.
Section \ref{sec:sep}). Obviously, if $z\in E$ and $w\in \DL
\subseteq L^1(\P)$ then $zw\in L^1(\P)$. Thus the bilinear form
$\ip\cdot\cdot:(E,\DL)\rightarrow\R$ defined by
\[ \ip z{w}:=\E{zw} \]
is well defined, making $(E,\DL)$ a bilinear system. From now on,
all polar wedges are defined with respect to this bilinear system.
Therefore, their structure depends on the choice of the wedge $K$
and the set $M$. Note that $E$ is a vector sublattice of $L^0$ which
contains the wedge $K$ (cf. \eqref{eqn:cone K in lattice}).

A standard measure theoretic argument shows that the linear
functionals $E\ni z\mapsto\E{zw}$ are non-degenerate for each $w\in
\DL\setminus\{0\}$.
Hence, the bilinear system $(E,\DL)$ is a left dual system.
Consequently, $\DL$ is the topological dual of $E$ under the weak
topology $\sigma(E,\DL)$:
\begin{equation}
\DL = (E, \sigma(E,\DL))'
\end{equation}
Note however that in general $(E, \DL)$ may not be a right dual
system (e.\,g., if $M$ consists of martingale measures only or if
$K\neq\{0\}$ and $\DL$ is finite dimensional), in which case $E$ is
not the dual of $\DL$ under the weak topology $\sigma(\DL,E)$. In
particular, the topology $\sigma(E,\DL)$ may not be Hausdorff.

Let $X\in E$ be an arbitrary contingent claim and $C \subseteq E$ be
an arbitrary wedge which contains the (fixed) wedge $K$. Recall that
\[
s_E(C) = \bigcap_{\Q \in M}(C - L^1_+(\Q)) = C - E_+ \subseteq E
\]
describes the umbrella hull of $C$ (cf. Proposition \ref{thm:rep of
umbrella cones in sublattices}). Consider the set
\[
A_X(C) : = \{x\in\R:X \leq x + G \textup{ for some } G \in C \}
\]
and assume that $A_X(K) \not= \emptyset$. Then $A_X(K)$ is bounded
from below, and we have
\[
-\infty < \sup_{\Q\in M}\E[\Q]{X} \leq \inf(A_X(K)) < +\infty.
\]
Obviously, $A_X(K) \subseteq A_X(C)$. If $A_X(C)$ is bounded from
below, we put
\begin{equation*}
\suprep(X;C) : = \inf(A_X(C)).
\end{equation*}
Note that $\suprep(X;C) \leq \suprep(X;K)$.
\begin{definition}\label{def:M based suprep price}
Let $M$ be a non-empty subset of $M_1 = M_1(\P; K)$ and $C \subseteq
E : = \bigcap_{\Q\in M}L^1(\Q)$ be a wedge which contains the wedge
$K$. Let $X \in E$. If $A_X(C)$ is non-empty and bounded from below
then the real number
\begin{equation*}
\suprep(X; C):= \inf(A_X(C)) = \inf\{x\in\R:X \leq x + G \textup{
for some } G \in C \}
\end{equation*}
is called $M$-based super-replication price of $X$.
\end{definition}
Due to Proposition \ref{thm:rep of umbrella cones in sublattices},
it follows immediately that
\begin{equation}\label{eqn:umbrella rep of M based suprep price}
\suprep(X;C) = \inf\big\{ x\in\R : X - x \in {s}_{E}(C)\big\} =
\inf\big\{ x\in\R : X - x \in \bigcap_{\Q \in M} \big(C -
L^1_+(\Q)\big)\big\}.
\end{equation}
Hence, $\suprep(X;C) = \suprep(X;{s}_{E}(C))$. If $M = M_\Phi$,
where $\Phi$ is the convex conjugate of a utility function, we call
$\suprep(X; C_\Phi)$ a utility-based super-replication price of $X$
(cf. \eqref{eqn:defcphi}). Note that $\suprep(X;C_{\textup{id}})$ is
the weak super-replication price of \cite{BiagFrit04}.
\comment{Thus, we have arrived at the following natural question:
Given an arbitrary wedge $K \subseteq L^0$, which wedges $C
\supseteq K$ and which non-empty sets $M$ of separating probability
measures would then be suitable to solve our pricing problem? Would
then $s_{E}(K)$ be a subset of $C$? }

The following important result reveals that the dual wedge of the
umbrella wedge $s_E(K)$ is contained in the positive wedge of
$L^1(\P)$ and can be represented with the help of suitable
\textit{probability} measures. It is not only crucial for our
analysis in this paper. By choosing the ``right'' dual pair namely,
it opens up a door to a canonical duality and a general
representation of $M$-based super-replication prices (cf.
\cite{OertelOwen07}).

\begin{lemma}\label{thm:rep of dual umbrella cones}
Let $K$ be an arbitrary wedge in ${L}^0$, and let $\emptyset \not= M
\subseteq M_1 = M_1(\P; K)$. Then
\[
\DL \cap \; \cone(R(M_1)) = s_E(K)^\lhd = (-L^\infty_+(\P))^\lhd
\cap K^\lhd = L_+^1(\P) \cap K^\lhd.
\]
In particular, the set $\DL \cap \cone(R(M_1))$ is $\sigma(\DL,
E)$-closed, and it contains $\overline{\cone(R(M))}^{\sigma(\DL,
E)}$.
\end{lemma}
\begin{proof}
First, we show that
\[
\DL \cap \; \cone(R(M_1)) \subseteq s_E(K)^\lhd \subseteq
(-L^\infty_+(\P))^\lhd \cap K^\lhd = L_+^1(\P) \cap K^\lhd.
\]
To this end, let $\Q \in M_1$ such that $R(\Q) = \displaystyle\radnik \in \DL$ 
and $X \in s_E(K)$. Then there exists $Y \in K \subseteq L^1(\Q)$
such that $X \leq Y$ and
\[
\langle X, \radnik \rangle = \E[\Q]{X} \leq \E[\Q]{Y} \leq 0.
\]
Thus, $\DL \cap R(M_1) \subseteq s_E(K)^\lhd$. Since $\emptyset
\not= R(M) \subseteq \DL \cap R(M_1)$, we may apply Lemma
\ref{thm:cone_of cone cap cx set} to the wedge $\DL$, and it
consequently follows that
\begin{equation}\label{eqn:towards the polar cone of an umbrella I}
\DL \cap \; \cone(R(M_1)) = \cone(\DL \cap R(M_1)) \subseteq
s_E(K)^\lhd.
\end{equation}
Since $L^\infty(\P) \subseteq L^1(\Q)$ for any $\Q \in M_1$, we have
$L^\infty(\P) \subseteq E$ and hence
\begin{equation}\label{eqn:towards the polar cone of an umbrella II}
-L^\infty_+(\P) \cup K \subseteq K - E_+ = s_E(K) \subseteq E.
\end{equation} Since
\begin{eqnarray*} (-L^\infty_+(\P))^\lhd & =
& \DL \cap \big\{W \in L^1(\P) :
\E[\P]{WV} \geq 0 \mbox{ for all } V \in L^\infty_+(\P) \big\}\\
& \subseteq & \DL \cap \big\{W \in L^1(\P) : \E[\P]{W\ind_{\{W <
0\}}} \geq 0 \big\}\\
& = & \DL \cap \big\{W \in L^1(\P) : \E[\P]{W^-} = 0 \big\}\\
& = & \DL \cap L^1_+(\P)\\
& \subseteq & (-L^\infty_+(\P))^\lhd,
\end{eqnarray*}
Lemma \ref{thm:permanence properties of dual cones} and inclusion
\eqref{eqn:towards the polar cone of an umbrella II} imply that
\begin{equation}\label{eqn:towards the polar cone of an umbrella III}
s_E(K)^\lhd \subseteq (-L^\infty_+(\P))^\lhd \cap K^\lhd = L^1_+(\P)
\cap K^\lhd.
\end{equation}
Hence,
\[
\DL \cap \; \cone(R(M_1)) \subseteq s_E(K)^\lhd \subseteq
(-L^\infty_+(\P))^\lhd \cap K^\lhd = L_+^1(\P) \cap K^\lhd.
\]
It remains to prove that $L^1_+(\P) \cap K^\lhd \subseteq \DL \cap
\; \cone(R(M_1))$. To this end, let $Y \in L^1_+(\P) \cap K^\lhd =
L^1_+(\P) \cap K^\lhd \cap \DL$.
Then $0 \leq Y = \sum_{k = 1}^n \alpha_k \radnik[\Q_k]$, for some $n
\in \N$, $\alpha_1, \ldots , \alpha_n \in \R$ and $\Q_1, \ldots ,
\Q_n \in M \subseteq M_1$, and $\E[\P]{XY} \leq 0$ for all $X \in
K$. Consider $\mu : = \sum_{k = 1}^n \alpha_k \Q_k$. Since each
$\Q_k$ in particular is a probability measure, it follows that
\[
0 \leq \E{Y} = \sum_{k = 1}^n \alpha_k = \mu(\Omega).
\]
If $\mu(\Omega) = \E{Y} = 0$ then $Y = 0 \in \DL \cap \;
\cone(R(M_1))$. If $\mu(\Omega) > 0$, it follows that $\Q : =
\sum_{k = 1}^n \frac{\alpha_k}{\mu(\Omega)} \Q_k$ is a probability
measure, $\Q \ll \P$ and
\[
\mu(\Omega)\,\E[\Q]{X} = \E[\P]{XY} \leq 0
\]
for all $X \in K$. Consequently, $\Q \in M_1$, and we obtain
\[
Y = \sum_{k = 1}^n \alpha_k \radnik[\Q_k] = \mu(\Omega) \radnik[\Q]
= \mu(\Omega) R(\Q) \in \DL \cap \; \cone(R(M_1)).
\]
\qed
\end{proof}
\begin{theorem}\label{thm:symmetric duality}
Let $K$ be an arbitrary wedge in ${L}^0$ and let $\emptyset \not= M
\subseteq M_1 = M_1(\P; K)$. Put $c_E(K) : =
\overline{s_E(K)}^{\sigma(E, \DL)}$. Then the following statements
are equivalent:
\begin{enumerate}
\item[$(i)$] $\overline{\cone(R(M))}^{\sigma(\DL, E)} = \DL \cap \cone(R(M_1))$;
\item[$(ii)$] $R(M)^\lhd = c_E(K)$;
\item[$(iii)$] $\overline{\cone(R(M))}^{\sigma(\DL, E)} = (c_E(K))^\lhd = s_E(K)^\lhd$.
\end{enumerate}
\end{theorem}
\begin{proof}
Assume that (i) holds. Then Proposition \ref{thm:bipolar_cones} and
Lemma \ref{thm:rep of dual umbrella cones} both together imply that
\begin{eqnarray*}
R(M)^\lhd & = & (R(M)^{\lhd\lhd})^\lhd =
\big(\overline{\cone(R(M))}^{\sigma(\DL,
E)}\big)^\lhd\\
& \stackrel{(i)}{=} & \big( \DL \cap \cone(R(M_1))\big)^\lhd =
(s_E(K)^\lhd)^\lhd = c_E(K),
\end{eqnarray*}
and (ii) follows.

Now assume that (ii) holds. Again, we may apply Proposition
\ref{thm:bipolar_cones} and Lemma \ref{thm:rep of dual umbrella
cones} and obtain
\[
\overline{\cone(R(M))}^{\sigma(\DL, E)} = R(M)^{\lhd\lhd}
\stackrel{(ii)}{=} c_E(K)^\lhd = (s_E(K)^{\lhd\lhd})^\lhd =
s_E(K)^\lhd.
\]
Hence, statement (iii) follows.

If (iii) is true, Lemma \ref{thm:rep of dual umbrella cones}
immediately implies (i).\qed
\end{proof}
A natural question is to ask for suitable candidates $M$ which
satisfy condition (ii) and hence the equivalent relations (i)
respectively (iii) of Theorem \ref{thm:symmetric duality}. In fact,
we shall recognise again that one candidate of this type is given by
$\widehat M_\Phi$ (cf. \cite{OertelOwen07}).
\begin{proposition}\label{thm:symmetric duality for faces}
Let $K$ be an arbitrary wedge in ${L}^0$, and let $\emptyset \not=
M$ be a face of $M_1 = M_1(\P; K)$. Then
\[
R(M)^\lhd = c_E(K).
\]
\end{proposition}
\begin{proof}
Firstly, due to Proposition \ref{thm:faces of M1 and their linear
hull} we have
\[
R(M) = R([M] \cap M_1) = \DL \cap R(M_1),
\]
so that
\[
\cone(R(M)) = \DL \cap \cone(R(M_1))
\]
(by Lemma \ref{thm:cone_of cone cap cx set}). Consequently, Lemma
\ref{thm:rep of dual umbrella cones} implies
\[
\cone(R(M))^\lhd = (s_E(K)^{\lhd})^\lhd = s_E(K)^{\lhd\lhd}.
\]
Now, we may apply Proposition \ref{thm: representation of a
polar_cone} and Proposition \ref{thm:bipolar_cones}, and it
consequently follows that
\[
R(M)^\lhd = \cone(R(M))^\lhd = s_E(K)^{\lhd\lhd} =
\overline{s_E(K)}^{\sigma(E,\DL)} = c_E(K).
\]
\qed
\end{proof}
\begin{corollary}\label{thm:general superreplication result} Let
$K$ be an arbitrary wedge in ${L}^0$, and let $\emptyset \not= M
\subseteq M_1 = M_1(\P, K)$ such that
\[
\overline{\cone(R(M))}^{\sigma(\DL, E)} = \DL \cap \cone(R(M_1)),
\]
respectively
\[
R(M)^\lhd = c_E(K).
\]
Let $X \in E$. Then the set $A_X(c_E(K))$ is bounded from below and
\[
\suprep(X; c_E(K)) = \sup_{\Q\in M}\E[\Q]X\,. 
\]
\end{corollary}
\begin{proof}
Since obviously
\[
c_E(K) - E_+ = \overline{s_E(K)}^{\sigma(E, \DL)} - E_+ \subseteq
\overline{K-E_+}^{\sigma(E, \DL)} = c_E(K),
\]
$c_E(K)$ is an umbrella wedge in E. Hence, Theorem
\ref{thm:symmetric duality} implies that
\begin{eqnarray*}
\suprep(X; c_E(K)) \!\!& = & \!\! \inf\{x \in \R : X - x \in c_E(K) \}\\
\!\!& = & \!\! \inf\{x \in \R : X - x \in R(M)^\lhd \}\\
\!\!& = & \!\! \inf\{x \in \R : \E[\Q]{X-x} \leq 0 \,\forall\,\Q \in M\}\\
\!\!& = & \!\! \inf\{x \in \R : \E[\Q]{X} \leq x \;\forall\; \Q \in M\}\\
\!\!& = & \!\! \sup_{\Q\in M}\E[\Q]X,
\end{eqnarray*}
and the claim follows.\qed
\end{proof}
\begin{corollary}
Let $K$ be an arbitrary wedge in ${L}^0$. If $\Phi$ satisfies the
growth condition \eqref{eqn:growth} and $\widehat M_\Phi \not=
\emptyset$ then
\[
R(\widehat M_\Phi)^\lhd = c_{\widehat{E_\Phi}}(K),
\]
where the polarisation $\lhd$ is based on the bilinear system
$(\widehat{E_\Phi}, \DLCW)$ which is given by $\widehat{E_\Phi} =
\bigcap_{\Q\in \widehat M_\Phi}L^1(\Q)$ and $\DLCW = [R(\widehat
M_\Phi)]$.
\end{corollary}
\begin{proof}
The result immediately follows from Lemma \ref{thm:likefrit} and
Proposition \ref{thm:symmetric duality for faces}, applied to $M : =
\widehat M_\Phi$ and the related bilinear system $(E, \DL) : =
(\widehat{E_\Phi}, \DLCW)$ (cf. \ref{eqn:bilinear system}).\qed
\end{proof}
Note that Lemma \ref{thm:rep of dual umbrella cones} and Theorem
\ref{thm:symmetric duality} hold for \textit{any} non-empty subset
$M$ of $M_1$. In particular, we may apply them to $M \in \{M_\Phi,
\widehat{M_\Phi}\}$ - presupposed that $M_\Phi \not= \emptyset$
respectively $\widehat{M_\Phi} \not= \emptyset$. The case $M =
\widehat{M_\Phi}$ is discussed in detail in \cite{OertelOwen07}. We
now apply our general results to the case $M = M_\Phi$ and provide a
well-known representation of the umbrella wedge $c_E(K)$ in Theorem
\ref{thm:weakclose}, implying a generalisation of Theorem 5 in
\cite{BiagFrit04}.

\section{The Special Case $M = M_\Phi$}\label{sec:M_Phi}
As before, we fix an arbitrary wedge $K$ in $L^0$. Assume that
$\Phi$ satisfies the growth condition \eqref{eqn:growth} and
$M_\Phi\neq\emptyset$. Let $E_\Phi := \bigcap_{\Q\in M_\Phi}L^1(\Q)$
be the vector space of all $M_\Phi$-integrable contingent claims.
Recall that $K \subseteq \bigcap_{\Q\in M_1}L^1(\Q) \subseteq
E_\Phi$. Consider the wedge
\[
K_\Phi : = {s}_{E_\Phi}(K) = K - (E_\Phi)_+  = \bigcap_{\Q\in
M_\Phi}\big(K-L^1_+(\Q)\big)
\]
of all $M_\Phi$-integrable contingent claims that can be dominated
by a terminal wealth in $K$. As shown in the previous section,
the linear span of the Radon-Nikodym derivatives of all probability
measures in $M_\Phi$ plays a fundamental role in the calculation of
super-replication prices of unbounded contingent claims in $E_\Phi$.
More precisely, we now consider the pair of vector spaces
\begin{equation}
  E_\Phi : = \bigcap_{\Q\in M_\Phi}L^1(\Q)\qquad\textup{and}\qquad
\DLC : = [R(M_\Phi)].
\end{equation}
It is well-known that the wedge $K^{\adm}$ of terminal wealths
arising from zero-financed admissible trading strategies is not
large enough for the purposes of a duality theory when considering
unbounded wealth. Similarly, in our general setting the wedge
$K_\Phi$ may not be large enough in order to obtain a dual
relationship of type \eqref{eqn:arbcone}. Along the lines of
equation \eqref{eqn:Cdefforadm} we therefore define the larger wedge
\begin{equation}\label{eqn:defcphi}
  C_\Phi:=\bigcap_{\Q\in M_\Phi}\normcl{K-L_+^1(\Q)},
\end{equation}
which in fact will turn out to be the closure of $K_\Phi$ with
respect to the weak topology $\sigma(E_\Phi, \DLC)$ (see Theorem
\ref{thm:weakclose}). Both, $K_\Phi$ and $C_\Phi$ are umbrella
wedges in $E_\Phi$, and
\begin{equation*}
K \subseteq K_\Phi \subseteq C_\Phi \subseteq E_\Phi.
\end{equation*}
We interpret $C_\Phi$ as the wedge of contingent claims which can be
approximated by super-replicable claims, where the investor has only
utility-induced restrictions on wealth.

The following result shows that the wedges $C_{\widehat\Phi}$ and
$C_\Phi$ are identical (c.f. \eqref{eqn:defphihat}).
\begin{lemma}\label{thm:redundant} If $\Phi$ satisfies the
growth condition \eqref{eqn:growth} and $M_\Phi\neq\emptyset$ then
$C_{\widehat\Phi} = C_\Phi$.
\end{lemma}

\begin{proof} Since $M_\Phi \subseteq\widehat M_\Phi =
M_{\widehat\Phi}$ it follows that $C_{\widehat\Phi} \subseteq
C_\Phi$. Therefore it suffices to show that for any $\Q_1\in\widehat
M_\Phi$ there exists $\Q\in M_\Phi$ such that
$\normcl{K-L_+^1(\Q)}\subseteq\normcl[\Q_1]{K-L_+^1(\Q_1)}$. Indeed,
let $\Q_1\in\widehat M_\Phi$, take any $\Q_0\in M_\Phi$ and define
$\Q=\frac12\Q_0+\frac12\Q_1$. Then $\Q\in M_\Phi$ (due to Lemma
\ref{thm:likefrit}). If $X\in\normcl{K-L_+^1(\Q)}$ then there exists
$\widetilde X_n\in K$ and $R_n\in L_+^1(\Q)$ such that
$X_n:=\widetilde X_n-R_n\overset{L^1(\Q)}\longrightarrow X$. Since
$\|R_n\|_{L^1(\Q_1)}\le2\|R_n\|_{L^1(\Q)}$ we have $X_n\in K -
L_+^1(\Q_1)$. Moreover,
$\|X-X_n\|_{L^1(\Q_1)}\le2\|X-X_n\|_{L^1(\Q)}\rightarrow0$ as
$n\rightarrow\infty$.\qed
\end{proof}

\begin{remark}\label{thm:largelosses}
Lemma \ref{thm:redundant} reveals an interesting economic insight
into the wedge $C_\Phi$. On inspection of definition
\eqref{eqn:defcphi}, one is lead to believe that $C_\Phi$ is highly
dependent on $\Phi$, and therefore on $U$. However, as a result of
Lemma \ref{thm:redundant}, we can replace in this definition the set
$M_\Phi$ by the set $M_{\widehat\Phi}=\widehat M_\Phi$ of pricing
measures with finite \textit{loss}-entropy.

The definition of the loss-entropy of a pricing measure only depends
upon the conjugate function $\Phi(y)$ for arbitrarily large values
of $y$ (see the discussion after equation \eqref{eqn:measfinitele}).
In turn, the behaviour of $\Phi(y)$ for large values of $y$
corresponds to the behaviour of the utility function $U(x)$ for
large negative values of $x$. Therefore, although the trader is
restricted in their choice of terminal wealths by their utility
function, this restriction actually depends only upon the investor's
preferences towards asymptotically large losses.
\end{remark}
\begin{theorem}\label{thm:polar2}
If $\Phi$ satisfies the growth condition \eqref{eqn:growth} and
$M_\Phi\neq\emptyset$ then
\begin{equation*}  C_\Phi = R\big(\widehat M_\Phi\big)^\lhd =
R(M_\Phi)^\lhd,
\end{equation*}
where the polarisation $\lhd$ is based on the bilinear system
$(E_\Phi, \DLC)$. In particular, $C_\Phi$ is $\sigma(E_\Phi,
\DLC)$-closed.
\end{theorem}
\begin{proof}
Let $\Q\in\widehat M_\Phi = M_{\widehat\Phi}$ and $X\in C_\Phi$.
Then $X\in C_{\widehat\Phi} \subseteq \normcl{K-L_+^1(\Q)}$. Thus
there exists a sequence $(X_n)_{n=1}^\infty\subseteq K-L_+^1(\Q)$
such that ${\| X - X_n \|}_{L^1(\Q)} \to 0$. Since $\Q\in M_1$, it
follows that $\E[\Q]{X_n}\leq 0$ for all $n$ and hence $\E[\Q]X\leq
0$, implying that
\begin{equation}\label{eqn:onedir} C_\Phi\subseteq R\big(\widehat
M_\Phi\big)^\lhd \subseteq R(M_\Phi)^\lhd. \end{equation} It remains
to show that $R(M_\Phi)^\lhd\subseteq C_\Phi$. To this end, we
proceed along the lines of the proof of the Kreps-Yan Theorem (cf.
\cite[Theorem 3.5.8]{ElliKopp}) and consider an arbitrary $Z\in
E_\Phi$ such that $Z\not\in C_\Phi$. Then there exists $\Q^\ast\in
M_\Phi$ such that $Z\not\in\normcl[\Q^\ast]{K - L_+^1(\Q^\ast)}$. By
the Hyperplane Separation Theorem, there exists a continuous linear
functional on $L^1(\Q^\ast)$ that separates $Z$ from the closed
wedge $\normcl[\Q^\ast]{K - L_+^1(\Q^\ast)}$, i.\,e., there exists
$\Lambda\in L^\infty(\Q^\ast)$ such that
\begin{equation*} \E[\Q^\ast]{\Lambda X}\le0<\E[\Q^\ast]{\Lambda Z} \end{equation*}
for all $X\in K - L_+^1(\Q^\ast)$. Since $-\ind_A\in K -
L_+^1(\Q^\ast)$ for any set $A \in {\mathcurl{F}}_T$, we deduce
that $\Lambda \geq 0$ $\Q^\ast$-a.\,s. and $\E[\Q^\ast]{\Lambda} >
0$. Thus, if we put $\Lambda^\ast = \Lambda
/\E[\Q^\ast]{\Lambda}$, then
\begin{equation*} \Q_0(A) := \E[\Q^\ast]{\Lambda^\ast\ind_A} \end{equation*}
defines a probability measure ${\Q}_0\ll\P$ on
$(\Omega,{\mathcurl{F}}_T)$, and the above inequality implies that
${\Q}_0\in M_1$ and $\E[{\Q}_0]{Z}>0$. Due to Lemma
\ref{thm:truncation_inequality} and the growth condition
\eqref{eqn:growth} of $\widehat\Phi$, it therefore follows that
\begin{equation*}
  \E{\widehat\Phi^+\left(\radnik[{\Q}_0]\right)} =
  \E{\widehat\Phi^+\left(\Lambda^\ast\radnik[\Q^\ast]\right)}
  \leq \widehat\Phi(0)^+ +
  \E{\widehat\Phi^+\left(\|\Lambda^\ast\|_{L^\infty(\Q^\ast)}\radnik[\Q^\ast]\right)}
  < \infty.
\end{equation*}
Hence, ${\Q}_0\in\widehat M_\Phi$.

Now pick any $\Q_1\in M_\Phi$. For $\varepsilon\in(0,1)$, Lemma
\ref{thm:likefrit} shows that the measure
$\Q_\varepsilon:=\varepsilon\Q_1+(1-\varepsilon)\Q_0\in M_\Phi$.
Moreover, if $\varepsilon$ is small enough then clearly
$\E[\Q_\varepsilon]Z>0$. Hence $Z \not \in R(M_\Phi)^\lhd$.\qed
\end{proof}
\begin{theorem}\label{thm:weakclose}
If $\Phi$ satisfies the growth condition \eqref{eqn:growth} and
$M_\Phi\neq\emptyset$ then
\begin{equation}\label{eqn:widehat Phi equals Phi dual cap M 1}
\DLC \cap \cone(R(M_1)) = \cone(R(\widehat M_\Phi)) =
\overline{\cone(R(M_\Phi))}^{\sigma(\DLC,E_\Phi)}.
\end{equation}
Moreover,
\begin{equation}\label{eqn:inL} C_\Phi =
c_{E_\Phi}(K) = \overline{K_\Phi}^{\sigma(E_\Phi,\DLC)}
\end{equation}
and
\begin{equation}
R(\widehat M_\Phi) = \overline{R(M_\Phi)}^{\sigma(\DLC,E_\Phi)}.
\end{equation}
\end{theorem}
\begin{proof}
Firstly, due to Proposition \ref{thm:newinput} we have
\[
R(\widehat M_\Phi) = R([M_\Phi] \cap M_1) = \DLC \cap R(M_1),
\]
so that
\[
\cone(R(\widehat M_\Phi)) = \DLC \cap \cone(R(M_1))
\]
(by Lemma \ref{thm:cone_of cone cap cx set}). Consequently, Lemma
\ref{thm:rep of dual umbrella cones} implies
\[
\cone(R(\widehat M_\Phi))^\lhd = (K_\Phi^{\lhd})^\lhd =
K_\Phi^{\lhd\lhd},
\]
where the polarisation $\lhd$ now is based on the bilinear system
$(E_\Phi, \DLC)$! Now, we may apply Theorem \ref{thm:polar2},
Proposition \ref{thm: representation of a polar_cone} and
Proposition \ref{thm:bipolar_cones}, and it follows
\[
C_\Phi = R\big(\widehat M_\Phi\big)^\lhd = \cone(R(\widehat
M_\Phi))^\lhd = K_\Phi^{\lhd\lhd} =
\overline{K_\Phi}^{\sigma(E_\Phi,\DLC)}.
\]
Consequently, Theorem \ref{thm:polar2} leads to
\begin{equation*} \cone(R(\widehat M_\Phi)) = \DLC \cap \cone(R(M_1))
= K_\Phi^{\lhd} = C_\Phi^\lhd = R(M_\Phi)^{\lhd\lhd} =
\overline{\cone(R(M_\Phi))}^{\sigma(\DLC,E_\Phi)}.
\end{equation*} Since $\ind_{\Omega} \in E_\Phi$, the functional
$\DLC\ni w\mapsto\E{w}$ is $\sigma(\DLC,E_\Phi)$-continuous,
implying that $\{w\in \DLC:\E{w}=1\}$ is
$\sigma(\DLC,E_\Phi)$-closed. Since
\begin{equation*} R(M_\Phi) = \cone(R(M_\Phi)) \cap \{w\in \DLC:\E{w}=1\},
\end{equation*}
we therefore obtain
\begin{align*}
  \overline{R(M_\Phi)}^{\sigma(\DLC,E_\Phi)}
  & \subseteq \overline{\cone(R(M_\Phi))}^{\sigma(\DLC,E_\Phi)}
  \cap\{w \in  \DLC:\E{w}=1\}\\
  & =   \cone(R(\widehat M_\Phi)) \cap \{w\in
  \DLC:\E{w}=1\}\\
  & =  \DLC \cap \cone(R(M_1)) \cap \{w\in
  \DLC:\E{w}=1\}\\
  & =  \DLC \cap R(M_1) = R(\widehat M_\Phi ).
\end{align*}
Remember that the last equality follows from Proposition
\ref{thm:newinput}. To prove the other inclusion, let $\Q_1 \in
\widehat M_\Phi$ arbitrary and fix $\Q_0 \in M_\Phi$. Then, due to
Lemma \ref{thm:likefrit}, $\Q^{(n)} : = (1-1/n) \Q_1 + 1/n \, \Q_0
\in M_\Phi$ for any $n \in \N$, and obviously $R(\Q^{(n)})
\longrightarrow R(\Q_1)$ in the topology $\sigma(\DLC,E_\Phi)$. \qed
\end{proof}
\begin{corollary}\label{thm:two}
Let $X \in E_\Phi$ and assume that $\Phi$ satisfies the growth
condition \eqref{eqn:growth} and $M_\Phi\neq\emptyset$. Then
\begin{equation}\label{eqn:inhatM}
  \suprep(X;C_\Phi) = \sup_{\Q\in M_\Phi}\E[\Q]X.
\end{equation}
\end{corollary}
\begin{proof}
Firstly, recall from equation \eqref{eqn:widehat Phi equals Phi dual
cap M 1} of Theorem \ref{thm:weakclose} that
\[
\DLC \cap \cone(R(M_1)) =
\overline{\cone(R(M_\Phi))}^{\sigma(\DLC,E_\Phi)}.
\]
Moreover, due to equation \eqref{eqn:inL} of Theorem
\ref{thm:weakclose}, we know that $C_\Phi = c_{E_\Phi}(K)$.
Therefore, we may apply Corollary \ref{thm:general superreplication
result} to $M = M_\Phi$. \qed
\end{proof}
\begin{remark}\label{thm:laterremark}
If $\Phi$ corresponds to a utility function which is supported on a
half-line (i.\,e., $U:(a,\infty)\rightarrow\R$, where $a>-\infty$)
then the convex conjugate $\Phi$ of $U$ is asymptotically linear as
$y\rightarrow\infty$. As mentioned in Remark
\ref{thm:pricingremarks} (ii), this means that $\widehat
M_\Phi=M_1$. For such utility functions therefore, the set of
super-replicable contingent claims does not depend specifically on
the shape of $U$. In fact, $C_\Phi=C_{id}$ where $id(y):=y$, and we
recover the polar relations
\begin{equation}
  R(M_1)^\lhd=C_{id}\qquad\textup{and}\qquad(C_{id})^\lhd = \cone(R(M_1)).
\end{equation}
This polarity is of a similar nature to \cite[Theorem
3.1]{KramScha99}, in the sense that it is {\it{utility
independent}}.
\end{remark}
\begin{remark}
Note that in equation \eqref{eqn:inL}, we may in fact take the
closure in any admissible topology (i.\,e., in any topology which is
stronger than the weak topology $\sigma(E_\Phi,\DLC)$ and weaker
than the Mackey topology $\tau(E_\Phi,\DLC)$). See \cite[\S{}98 and
\S{}103]{Heus82} for an explanation of the details.
\end{remark}

\section{The Case of Admissible Trading
Strategies}\label{sec:admiss case}

In this section we consider the particular case where $K=K^{\adm}$
is the wedge of attainable terminal wealths resulting from
zero-financed {\it admissible} trading strategies. As an application
of our general framework, we show that in this case every contingent
claim in $C_\Phi = : C^{\adm}_\Phi$ even can be approximated by {\it
bounded} contingent claims which are dominated by terminal wealths
in $K^{\adm}$. This approximation is given with respect to the
($K^{\adm}$-related) weak topology $\sigma(E_\Phi,\DLC)$.

By
\begin{align*}\label{eqn:C}
C^{\adm} &:= (K^{\adm}-L_+^0)\cap L^\infty(\P)\\
  &=  \{X\in L^\infty(\P): X\le \widetilde X \textup{ for
  some } \widetilde X \in K^{\adm}\}
\end{align*}
we denote the wedge of all a.\,s {\it bounded} contingent claims
that can be dominated by a terminal wealth in $K^{\adm}$. Since
$K^{\adm} \subseteq K^{\adm}_\Phi$, it follows that
\[
C^{\adm} \subseteq (K^{\adm} - L_+^0) \cap L^1(\Q) \subseteq
K^{\adm}-L^1_+(\Q)
\]
for any $\Q \in M_\Phi$ and thus
\begin{equation*}
C^{\adm} \subseteq K^{\adm}_\Phi \subseteq C^{\adm}_\Phi \subseteq
E_\Phi.
\end{equation*}
A benefit of using admissible strategies is the following
approximation result:
\begin{lemma}\label{thm:convergence_lemma}
Let $X \in K^{\adm}$. Then there exists a constant $c \geq 0$ and a
sequence $(X_n)_{n \in \N}$ of random variables such that for any
probability measure $\Q \ll \P$ the following properties hold:
\begin{enumerate}
\item $X_n \in C^{\adm}$ for all $n \in \N$; \item $-c \leq X_1 \leq X_2
\leq \ldots \leq X \;\;\; {\Q}$-a.\,s; \item $X_n \uparrow X
\;\;\; {\Q}$-a.\,s; \item $-c \leq \lim\limits_{n \to
\infty}\E[\Q]{X_n} = \E[\Q]{X} \in \R \cup \{\infty\}$.
\end{enumerate}
\end{lemma}
\begin{proof}
Let $n \in \N$ and $\Q\ll\P$. Put $X_n : = X \wedge n$. Since $X \in
K^{\adm}$, there exists a real constant $c \geq 0$ such that $-c
\leq X_n \leq X$ $\P$-a.\,s, implying that $\Q$-a.\,s $(X_n)$ is a
monotonically increasing sequence in $C^{\adm}$ with limit $X$.
Hence, $0 \leq X_n + c \uparrow X +c$ $\Q$-a.\,s, and an application
of Lebesgue's Monotone Convergence Theorem finishes the proof. \qed
\end{proof}
As a consequence of Lemma \ref{thm:convergence_lemma} we may
substitute the wedge $K^{\adm}$ by the wedge $C^{\adm} \subseteq
L^\infty(\P)$ in equation \eqref{eqn:altdefM1} to get
\begin{equation}\label{eqn:replbyC}
R(M_1^\adm) = \{Y\in L_+^1(\P):\E Y=1\text{ and }\E{XY}\le0\textup{
for all }X\in C^{\adm}\}.
\end{equation}
\begin{theorem}\label{thm:admiss weakclose}
If $\Phi$ satisfies the growth condition \eqref{eqn:growth} and
$M_\Phi^\adm \neq\emptyset$ then
\begin{equation}\label{eqn:admiss inL}
C^{\adm}_\Phi = \overline {C^{\adm}}^{\sigma(E_\Phi,\DLC)}.
\end{equation}
\end{theorem}
\begin{proof}
Due to equation \eqref{eqn:replbyC} and Proposition
\ref{thm:newinput},
 $ {(C^\adm)}^\lhd \subseteq \DLC\cap\cone(R(M_1^\adm))=\cone(R(\widehat{M_\Phi^\adm}))$.
Hence, as a result of Theorem \ref{thm:polar2} and Proposition
\ref{thm:bipolar_cones},
\begin{equation*}
  \overline {C^{\adm}}^{\sigma(E_\Phi,\DLC)} \subseteq C^{\adm}_\Phi =
  (\cone(R(\widehat{M_\Phi^\adm})))^\lhd \subseteq {(C^\adm)}^{\lhd\lhd} = \overline
  {C^{\adm}}^{\sigma(E_\Phi,\DLC)}.\tag*{\qed}
\end{equation*}
\end{proof}

\section{Appendix}\label{sec:appendix}
Let us recall an important version of the Hyperplane Separation
Theorem in finite-dimensional vector spaces which is not only known
as one of the main building blocks for duality theorems in linear
programming. It also has other numerous applications, e.\,g., to the
Karush-Kuhn-Tucker theorem in nonlinear programming and zero-sum
games in economic theory (cf. \cite{BorLew06}).
\begin{theorem}[Farkas' Lemma]\label{thm:Farkas}
Let $m, n \in \N, A \in M(m \times n; {\R}^n)$ and $b \in {\R}^m$.
Then either
\begin{itemize}
\item[$($i$)$] there exists $x \in {\R}^n$ such that $x \geq 0$
and $Ax = b$
\end{itemize}
or
\begin{itemize}
\item[$($ii$)$] there exists $y \in {\R}^m$ such that $\langle y,
b \rangle > 0$ and $A^{\top}y \leq 0$.
\end{itemize}
\end{theorem}
Theorem~\ref{thm:Farkas} implies the non-trivial fact that the
finitely generated wedge $A({\R}^n_{+})$ is {\it closed} in
${\R}^m$. More generally, the following statement, which transfers
Theorem~\ref{thm:Farkas} to infinite-dimensional Banach lattices,
shows that Farkas' Lemma even is {\it equivalent} to the closedness
of the wedge $A({\R}^n_{+})$! Concerning the basics of Banach
lattices and positive operators, we refer the reader to
\cite{Schaef74}.
\begin{theorem}\label{thm:Farkas_infty}
Let $E, F$ be arbitrary Banach lattices and $A : E \longrightarrow
F$ be a continuous linear operator. Then the following statements
are equivalent:
\begin{itemize}
\item[$($i$)$] $A(E_+)$ is $\sigma(F, F')$-closed;
\item[$($ii$)$]
Let $b \in F$. Then either
\begin{itemize}
\item[$($1$)$] there exists $x \in E_+$ such that $Ax = b$
\end{itemize}
or
\begin{itemize}
\item[$($2$)$] there exists $y^\prime \in F^\prime$ such that
$\langle b, y^\prime \rangle > 0$ and $A^\prime y^\prime \leq 0$.
\end{itemize}
\end{itemize}
\end{theorem}
\begin{proof}
First assume that $(i)$ holds. We want to show $(ii)$. To this end,
we consider the wedge
\[
G : = ((A^\prime)^{-1}(-E^\prime_+))^\lhd  = \big\{b \in F : \langle
b, y^\prime \rangle \leq 0 \textup{ for all } y^\prime \in
(A^\prime)^{-1}(-E^\prime_+)\big\} \subseteq F.
\]
Since $F\setminus G$ exactly consists of all $b \in F$, satisfying
condition $(2)$ in $(ii)$ above, we only have to show that
\begin{equation}\label{eqn:Farkas_cone_1}
A(E_+) = G.
\end{equation}
Obviously, by definition of $G$, $A(E_+) \subseteq G$. Let $y^\prime
\in (A(E_+))^\lhd$. Then $\langle x, A^\prime y^\prime \rangle =
\langle Ax, y^\prime \rangle \leq 0$ for any $x \in E_+$. Hence,
$y^\prime \in (A^\prime)^{-1}(-E^\prime_+)$. Thus, $(A(E_+))^\lhd
\subseteq (A^\prime)^{-1}(-E^\prime_+)$ and consequently,
\begin{equation}\label{eqn:Farkas_cone_2}
A(E_+) \subseteq G \subseteq (A(E_+))^{\lhd\lhd} .
\end{equation}
Due to Proposition \ref{thm:bipolar_cones} and the assumed weak
closedness of the wedge $A(E_+)$, it follows that
\[
(A(E_+))^{\lhd\lhd} = \overline{\cone(A(E_+))}^{\sigma(F, F')} =
\overline{A(E_+)}^{\sigma(F, F')} = A(E_+)\,,
\]
Hence, we obtain equality \eqref{eqn:Farkas_cone_1}.

Assume now that statement $(ii)$ is true. We want to show $(i)$. To
this end, let $b \in \overline{A(E_+)}^{\sigma(F, F')}$. Then there
exists a net $(x_\alpha) \subseteq E_+$ such that $b$ is the
$\sigma(F, F')$-limit of the net $(Ax_\alpha)$. Assume by
contradiction that $b \notin A(E_+)$. Then, due to assumption
$(ii)$, condition $(2)$ must be true (since $(1)$ is false).
Consequently, there exists a continuous linear functional $y^\prime
\in F^\prime$ such that for any $x \in E_+$ we have
\[
\langle Ax, y^\prime \rangle = \langle x, A^\prime y^\prime
\rangle\leq 0 < \langle b, y^\prime \rangle.
\]
In particular, $\langle b, y^\prime \rangle = \lim_\alpha \langle
Ax_\alpha, y^\prime \rangle \leq 0$ - a contradiction.  Hence, $b
\in A(E_+)$\,.\qed
\end{proof}


\comment{
\begin{theorem}\label{thm:one}
If $\Phi$ satisfies the growth condition \eqref{eqn:growth}, and
$M_\Phi\neq\emptyset$ then
\begin{equation*} (C_\Phi)^\lhd = (K_\Phi)^\lhd =
L^1_+(\P) \cap K^\lhd = \cone(\widehat M_\Phi). \end{equation*} In
particular, $\cone(\widehat M_\Phi)$ is $\sigma(\DLC,
E_{\Phi})$-closed.
\end{theorem}
\begin{proof}
By polarising equation \eqref{eqn:onedir}, we see that
\begin{equation*} \cone(\widehat M_\Phi) \subseteq (\widehat
M_\Phi)^{\lhd\lhd}\subseteq(C_\Phi)^\lhd \subseteq (K_\Phi)^\lhd.
\end{equation*} Since $(-L_+^\infty(\P)) \cup K \subseteq
K_\Phi$, it follows from the definition of $M_1$ that
\begin{equation*} K_\Phi^\lhd \subseteq (-L_+^\infty(\P))^\lhd
\cap K^\lhd \subseteq L_+^1(\P)\cap K^\lhd =
\DLC\cap\cone(M_1)=\cone(\widehat M_\Phi),
\end{equation*} the last equality being the result of Proposition
\ref{thm:newinput}.\qed
\end{proof}
}


\end{document}